\documentstyle[psfig,11pt,amssymb]{article}

\begin{document}

\title{Option pricing with log-stable L\'{e}vy processes}

\author{Przemys{\l}aw Repetowicz and Peter Richmond
\footnote{Department of Physics, Trinity College Dublin 2, Ireland.}}


\date{}

\maketitle


\begin{abstract}
We model the logarithm of the price (log-price) of a financial asset as a random variable obtained by projecting an operator stable random vector with a scaling index matrix $\underline{\underline{E}}$ onto a non-random vector. The scaling index $\underline{\underline{E}}$ models prices of the individual financial assets (stocks, mutual funds, etc.). We find the functional form of the characteristic function of real powers of the price returns and we compute the expectation value of these real powers and we speculate on the utility of these results for statistical inference. Finally we consider a portfolio composed of an asset and an option on that asset. We derive the characteristic function of the deviation of the portfolio, \mbox{${\mathfrak D}_t^{({\mathfrak t})}$}, defined as a temporal change of the portfolio diminished by the the compound interest earned. We derive pseudo-differential equations for the option as a function of the log-stock-price and time and we find exact closed-form solutions to that equation. These results were not known before. Finally we discuss how our solutions correspond to other approximate results known from literature,in particular to the well known Black \& Scholes equation. 
\end{abstract}

\noindent{\small{\bf Key words and phrases:}\ Option pricing, heavy tails
, operator stable, fractional calculus.

\baselineskip=15pt

\section{Introduction}
Early statistical models of financial markets assumed 
that asset price returns are independent, identically distributed (iid) Gaussian variables. 
\cite{Bach}. 
However, evidence has been found \cite{GopiMeyer} that the returns
exhibit power law (fat) tails in the high end of the distribution.
Except at very high frequencies or short times (\cite{GopiMeyer}),  
a better statistical description for many financial assets 
is provided by a model where the logarithm of the price
is a heavy tailed one-dimensional L\'{e}vy $\mu$-stable process \cite{Mandel,Fama, portfolio, MittnikR}.
Since the tail parameter $\mu$ that measures the probability of large 
price jumps will vary from one financial asset to the next, 
a model based on operator stable L\'{e}vy processes \cite{Meerschaert} is appropriate. 
This model allows the tail index to differ for each financial asset in the portfolio.  
Hence we formulate a model where
the log-price is a projections of an operator stable random vector
onto a predefined
direction (this projection determines the portfolio mix).
The cumulative probability distribution of the log-price diminishes
as a mixture of power laws and 
thus the higher-order moments of the distribution may not exist 
and the characteristic function of the distribution may not be analytic. 

Due to the constraints on the size of this paper we only include new results
leaving proofs for further publications.


\section{The model of the stock market\label{sec:MarketModel}}
In this section we define the model.
In the following we recall certain known properties of operator stable distributions 
and we derive Fourier transforms of real powers of projections of 
operator stable vectors onto a non-random vector.
In subsections (\ref{sec:Scal}) and (\ref{sec:ScalRot}) we derive 
Fourier transforms of operator stable random vectors for 
particular forms of parameters of the distribution.

\subsection{The basic properties and new results\label{sec:TheStockPrice}}  
Let $\log(S_t)$ be the logarithm of the price of the portfolio (log-price) at time $t$. 
We assume that the temporal change of the log-price
is composed of two terms, a deterministic term and a fluctuation term viz:
\begin{equation}
\frac{dS_t}{S_t} = \frac{S_{t+dt} - S_t}{S_t} 
= \alpha d t + \vec{\sigma} \cdot d\vec{L}_{t}
\label{eq:StockPriceFluct}
\end{equation}
The parameters  $\alpha \in {\mathbb R}$ (the drift) 
and the elements of the $D$ dimensional vector
\mbox{$\vec{\sigma} := \left(\sigma_1,\dots,\sigma_D\right)$} 
(the portfolio mix) are assumed to be non-random constants.
The random vector $\vec{ L}_t$ is (strictly) operator stable, meaning that
it is an operator-normalized limit of a sum of some independent, identically
distributed (iid) random vectors $\vec{X}_i$. 
We have
\begin{equation}
\vec{L}_t:= \lim_{n \to\infty}
n^{-E}
\sum_{i=1}^{\lfloor n t \rfloor} \vec{X}_i
\label{eq:Fluctuation}
\end{equation}
where $n^{-E}=\exp(-E\ln n)$ and $E$ is a real $D$-dimensional matrix
such that the equality holds in distribution.
The class of distributions of the former vectors related to a given matrix $E$
is termed an attraction domain of an operator stable law.
Members of such class are usually unknown.


We now recall some known facts \cite{Meerschaert} concerning the 
operator stable probability
density \mbox{$\omega_{\vec{L}_t}$}  and its Fourier transform
\mbox{$\tilde{\omega}(\vec{k}) := {\mathcal F}_{\vec{x}}[\omega](\vec{k})$}.

The following identities hold:

\begin{equation}
\omega_{\vec{L}_t}(\vec{x})=
\omega_{\vec{L}_1}(t^{-E}\vec{x})\ \mbox{det$(t^{-E})$}\quad {\rm for}\ {\rm all}\ t>0. 
\label{eq:ISExpDef}
\end{equation}
\vskip10pt
and
\begin{equation}
\tilde{\omega}^t(\vec{k}) =
\exp(-t \phi(\vec{k}) ) =
\tilde{\omega}( t^{E^{T}} \vec{k} ) =
\tilde{\omega}_{\vec{L}_t}(\vec{k})
\label{eq:ScalingLimitCorFourTransf}
\end{equation}
where $E^{T}$ is the transpose of
$E$ and $\phi$ is the 
negative logarithmic characteristic function of the random vector
$\vec{L}_1$.
In the following we assume that that function is even:
\begin{equation}
\phi(\vec{k}) = \phi(-\vec{k})
\label{eq:EvenLogCharFct}
\end{equation}
The identities (\ref{eq:ISExpDef}) are termed as a self-similar property 
of the random walk $\vec{L}_t$.


A motivation for introducing model
(\ref{eq:StockPriceFluct}) is statistical inference of parameters
of a distribution that describes real financial data.
In this context it is useful to know analytically the distribution of a real power 
of the integrated fluctuation term
in (\ref{eq:StockPriceFluct}).  
In general this is not known.
Here we derive some new results for operator stable L\'{e}vy distributions.
Denote by $\nu^{(\beta)}_{t}(z)$ the pdf 
of a random process 
\mbox{$\Xi^{(\beta)}_{t} :=
\left(\vec{\sigma} \cdot \vec{L}_{t} \right)^\beta $}
and by
\mbox{$\tilde{\nu}^{(\beta)}_t(k)$} its Fourier transform.
For $\beta \ge 1$ the identity holds:

\begin{eqnarray}
\tilde{\nu}^{(\beta)}_t(k) =
\left\{
\begin{array}{rr}
\int\limits_{-\infty}^\infty d\lambda \,
\tilde{\omega}(t^{E^{T}}
\hat{\sigma} \lambda)
{\mathcal K}^{(\beta)}({\mathfrak k}, \lambda)
& \mbox{for $\beta \ne 1$} \\
\tilde{\omega}(t^{E^{T}} \vec{\sigma} k)
& \mbox{for $\beta = 1$}
\end{array}
\right.
\label{eq:FourTrNMarginalProp}
\end{eqnarray}
where
\begin{eqnarray}
 {\mathcal K}^{(\beta)}(k, \lambda) :=
\left\{
\begin{array}{rr}
\frac{1}{2\pi} \int\limits_0^\infty dt e^{-k t^\beta}
 \left(
  {\mathfrak r} e^{-\imath \lambda {\mathfrak r} t} + {\mathfrak r}\mu e^
{\imath \lambda {\mathfrak r} \mu t}
\right)
&
\quad\mbox{if $\lfloor \beta \rfloor$ is even}
\\
\frac{1}{2\pi} \int\limits_0^\infty dt e^{-k t^\beta}
 \left(
  {\mathfrak r} e^{-\imath \lambda {\mathfrak r} t} + \bar{{\mathfrak r}}
\mu e^{\imath \lambda \bar{{\mathfrak r}} \mu t}
\right)
&
\quad\mbox{if $\lfloor \beta \rfloor$ is odd}
\end{array}
\right.
\label{eq:KernelDefIb}
\end{eqnarray}
  and $\hat{\sigma} := \vec{\sigma}/\left|\vec{\sigma}\right|$
and
\mbox{${\mathfrak r} := \exp(\imath \pi/(2 \beta))$},
\mbox{$\mu := \exp(-\imath \pi \{ \beta \}/(\beta) )$}.
The symbols $\lfloor\beta \rfloor$ and $\{\beta \}$
mean the biggest integer not larger then $\beta$
and the fractional part of $\beta$ respectively,
and ${\mathfrak k} := k \sigma^\beta$.

In addition for even values of $\lfloor \beta \rfloor$ we have:
 \begin{equation}
\tilde{\nu}^{(\beta)}_t(k) =
\frac{1}{\pi}
\int\limits_{l({\mathfrak r}, \beta)}
dz e^{-z} 
\int\limits_0^\infty d\xi 
\left(\frac{\sin(\xi)}{\xi}\right)
\tilde{\omega}(t^{E^T} 
\left(\frac{k}{z}\right)^{1/\beta} \frac{\xi}{\mathfrak r}
\vec{\sigma}
)
\label{eq:FractMomCharacTFct}
\end{equation}
where the integration line $l({\mathfrak r}, \beta))$ reads:
\begin{equation}
l({\mathfrak r}, \beta) :=  
[0, {\mathfrak r}^{-\beta} \infty] 
\cup
[0, (-1)^\beta {\mathfrak r}^{-\beta} \infty] 
\label{eq:Integrationline}
\end{equation}
The identities (\ref{eq:FractMomCharacTFct}) and (\ref{eq:Integrationline})
may be useful for describing the magnitude of the fluctuations of a random walk.

Thus it follows that the fractional moments of the scalar product 
are obtained by differentiating the Fourier transform at $k=0$.
We will obtain closed form results for these moments 
in section (\ref{sec:ExistenceofMoments}).
Here we only recall that in the non-Gaussian case, due to (\ref{eq:EvenLogCharFct}), 
we have:
\begin{equation}
E\left[ (\vec{\sigma} \cdot \vec{L_t})^\beta \right] = 
\left\{
\begin{array}{cc} 
\infty & \mbox{when $\beta$ is even} \\ 
    0  & \mbox{when $\beta$ is odd}
\end{array}
\right\}
\label{eq:InfiniteMoments}
\end{equation}
In addition for $\beta>0$ there the moment exists only
if $\beta$ does not exceed a certain threshold value.

\vskip10pt

It is our objective to 
price options on the portfolio of stocks driven 
by operator stable fluctuations (see section (\ref{sec:OptionPrice})). 
By this we mean a theory that 1) 
allows inference of the the stable index $E$ and the L\'{e}vy measure of the whole vector of stock prices (market)from a statistical sample and 2)hedges against risk in the market by the construction of an appropriate option. To the best of our knowledge, this has not yet been achieved.

The generic properties of operator stable probability distributions
and their marginals are described in \cite{Meerschaert}.
Here we recall some known facts and we analyse two particular cases 
of the stable index.
The Fourier transform $\tilde{\omega}(\vec{k})$ 
is uniquely determined via the stable index $E$ and the 
log-characteristic function $\phi$ confined to a unit sphere.
This can be seen by representing the vector $\vec{k}$ in 
the Jurek coordinates viz
\begin{equation}
\vec{k} = r_{\vec{k}}^{E^{T}} \cdot \vec{\theta}_{\vec{k}}
\label{eq:JurekCoordinates}
\end{equation}
where
$\left|\vec{\theta}_{\vec{k}}\right| = 1$.
Using the scaling relation (\ref{eq:ScalingLimitCorFourTransf}) we get:
\begin{equation}
\exp(-r_{\vec{k}} \phi(\vec{\theta}_{\vec{k}}) ) = 
\tilde{\omega}^{r_{\vec{k}}}(\vec{\theta}_{\vec{k}}) = 
\tilde{\omega}(r_{\vec{k}}^{E^{(T)}} \cdot \vec{\theta}_{\vec{k}}) = 
\tilde{\omega}(\vec{k}) =
\exp(-\phi(\vec{k}) )
\label{eq:ScalRels}
\end{equation}
and thus
\begin{equation}
\phi(\vec{k}) = r_{\vec{k}} \phi(\vec{\theta}_{\vec{k}})
\label{eq:LogCharacteristicFctDef}
\end{equation}

Since it follows from the Jordan decomposition theorem that every matrix $E$ 
is, in a certain basis, a block diagonal matrix 
the set of all possible jump intensities $\omega$ is narrowed down 
to few classes of solutions only, each one corresponding to a particular 
Jordan decomposition of the matrix $E$.
We now firstly investigate a few classes of solutions
as a function of $E$ and subsequently the generic solution for an arbitrary $E$.
The existence results in this field are given in \cite{Meerschaert}.
We stress that, contrary to \cite{Meerschaert},
we aim at computing the characteristic
functions and the fractional moments in closed form
rather than only showing 
their existence.

\subsection{Pure scaling\label{sec:Scal}}
In this case $E = (D\mu)^{-1} I$ where $I$ is a
$D$ dimensional identity matrix and $\mu > 0$ is a constant.
From (\ref{eq:JurekCoordinates}) we see that:
\begin{equation}
\vec{k} = r_{\vec{k}}^{1/(D\mu)} \vec{\theta}_{\vec{k}}
= \left|\vec{k}\right| \frac{\vec{k}}{\left|\vec{k}\right|}
\label{eq:PureScal}
\end{equation}
hence $r_{\vec{k}} = \left|\vec{k}\right|^{D \mu}$
and $\vec{\theta} = \vec{k}/\left|\vec{k}\right|$ and so
\begin{equation}
\tilde{\omega}(\vec{k}) = 
\exp\left(-\left|\vec{k}\right|^{D \mu} \phi(\frac{\vec{k}}{\left|\vec{k}\right|})\right)
\label{eq:JumpPdfPureScaling}
\end{equation}

The $\beta$-marginal probability density function 
from (\ref{eq:FourTrNMarginalProp}) reads:
\begin{eqnarray}
\tilde{\nu}_t^{(\beta)}(k) =
\int_{-\infty}^\infty d\lambda
\exp\left\{
-  t |\lambda|^{D \mu}
\phi(\frac{\vec{\sigma}}{|\vec{\sigma}|}\mbox{sign$(\lambda)$})
\right\}
{\mathcal K}^{(\beta)}(k \sigma^\beta, \lambda)
\label{eq:PureScalnubeta}
\end{eqnarray}
where the kernel ${\mathcal K}$ is defined in
(\ref{eq:KernelDefIb}).
From the properties of the Gamma function 
we obtain easily the 
fractional moment of the scalar product as:
\begin{equation}
E\left[(\vec{\sigma}\cdot \vec{L}_t)^\beta\right] = 
{\mathfrak C}(\beta,D \mu)
\cdot
\left(\sigma t^{(D \mu)^{-1}}\right)^\beta
\cdot
\left(
\cos\left(\frac{\pi \beta}{2}\right) 
e^{\imath \frac{\beta \pi}{2} }
\right)
\cdot
\phi_\pm^{\beta/(D \mu)} 
\label{eq:FinalMomentResultPureScal}
\end{equation}
where
\begin{equation}
{\mathfrak C}(\beta,D \mu) :=
\left(
\frac{2^{\beta+1}}{D\mu \sqrt{\pi}}
\Gamma(\frac{\beta+1}{2})
\frac{\Gamma(-\frac{\beta}{D\mu})}{\Gamma(-\frac{\beta}{2})}
\right)
\label{eq:PreFactorFracMom}
\end{equation}
Here 
$\phi_\pm := \phi(\pm)$.
The moment exists for $\beta < D\mu$.
The prefactor (\ref{eq:PreFactorFracMom}) in (\ref{eq:FinalMomentResultPureScal})
fits in with the known result for the fractional moment of a modulus of a stable variable
(see equation (3.6) page 32 in \cite{NikiasShao}).
For the derivations of that result  by means of the Mellin-Stieljes transform see
\cite{Zolotarev,CambanisMiller}
and by means of characteristic functions see \cite{Wolfe}.

\subsection{Scaling \& rotation\label{sec:ScalRot}}
In this case $D=2$ and we chose:
\begin{equation}
E =
\left(
\begin{array}{cc}
(2\mu)^{-1} & -b          \\
b           & (2\mu)^{-1}
\end{array}
\right)
\end{equation}
Clearly the trace $\mbox{Tr$\left[E\right] = \mu^{-1}$}$.
We denote by
${{ O}}_{\beta} :=
\left(
\begin{array}{rr} \cos(\beta) & -\sin(\beta) \\
                  \sin(\beta) & \cos(\beta)
\end{array}
\right)$ a two dimensional rotation by an angle \mbox{$\beta \in {\mathbb R}$}.
The mapping:
\begin{equation}
r^{E^{T}} :
{\mathbb R}^2 \ni \vec{k}
\rightarrow
r^{(2 \mu)^{-1}}
{{ O}}_{-b \log(r)} \vec{k}
\in {\mathbb R}^2
\end{equation}
changes the length of $\vec{k}$ by a multiplicative factor 
$r^{(2 \mu)^{-1}}$
and rotates the vector by an angle $-b \log(r)$.
The Jurek coordinates read: 
$\vec{r}_{\vec{k}} = \left|\vec{k}\right|^{2 \mu}$
and 
\mbox{$\vec{\theta}_{\vec{k}} = O_{2 b \mu \log(\left|k\right|)} 
(\vec{k}/\left|\vec{k}\right|)$} and so
\begin{equation}
\tilde{\omega}(\vec{k}) = \exp\left(-\left|\vec{k}\right|^{2 \mu} 
\phi\left(O_{2 b \mu \log(\left|k\right|)} 
\frac{\vec{k}}{\left|\vec{k}\right|}\right)\right)
\label{eq:JumpPdfScalingRot}
\end{equation}

The $\beta$-marginal probability density function and the fractional moment read:
\begin{eqnarray}
\tilde{\nu}_t^{(\beta)}(k) =
\int_{-\infty}^\infty d\lambda
\exp\left\{
-  t |\lambda|^{2 \mu}
\phi(
O_{2 b \mu \log(|\lambda|)}
\frac{\vec{\sigma}}{|\vec{\sigma}|}
)
\right\}
{\mathcal K}^{(\beta)}(k \sigma^\beta, \lambda)
\label{eq:CharacterFctScalRot}
\end{eqnarray}
and
\begin{equation}
E\left[\left(\vec{\sigma}\cdot \vec{L}_t\right)^\beta\right] =
{\mathfrak C}(\beta,2\mu)
\cdot
\left(
\sigma t^{(2 \mu)^{-1}}
\right)^\beta
\cdot
\left(
\cos\left(\frac{\pi \beta}{2}\right) 
e^{\imath \frac{\beta \pi}{2} }
\right)
\cdot
\frac{1}{2\pi}
\int\limits_0^{2\pi} d\eta \phi^{\beta/(2\mu)}	(O_\eta \hat{\sigma})
\cdot
\label{eq:FracMomScalRot}
\end{equation}
respectively.
Here 
${\mathfrak C}(\beta,2\mu)$ is defined in (\ref{eq:PreFactorFracMom}).
The moment exists for $\beta < 2\mu$.
Since the moment depends on the average of a power of the log-characteristic function
over the unit sphere we
conclude that the knowledge of the moments
does not determine the distribution in a unique manner.

In the following section we compute the
fractional moments of the scalar product
in the generic case of a operator stable distribution.

%
%
%
\subsection{The generic case \label{sec:ExistenceofMoments}}
Assume that the stable index has $D$ different eigenvalues 
$\left\{\lambda_p\right\}_{p=1}^D$
that are either real or pairwise complex conjugate. 
Then the following spectral decomposition holds:
\begin{equation}
\underline{\underline{E}}^T =
\underline{\underline{O}} \cdot 
\mbox{Diag$\left(\lambda\right)$} \cdot \underline{\underline{O}}^{-1}
\label{eq:SpectralDecompStablIndex}
\end{equation}
where 
\mbox{$\mbox{Diag$\left(\lambda\right)$} 
:= \left( \delta_{i,j} \lambda_j \right)_{i,j=1}^D$}
and such that the matrix $\underline{\underline{O}}$ is unitary
\begin{equation}
\underline{\underline{O}} \cdot \underline{\underline{O}}^{\dagger} = 
\underline{\underline{O}}^{\dagger} \cdot \underline{\underline{O}} = 1
\end{equation}
Then from (\ref{eq:SpectralDecompStablIndex}),
from the definition of the operator 
$r^{\underline{\underline{E}}^T} := \exp(\underline{\underline{E}}^T \log(r))$
and from the Cayley-Hamilton theorem we easily arrive at the identity:
\begin{equation}
r^{\underline{\underline{E}}^T} :=
\underline{\underline{O}} \cdot 
\mbox{Diag$\left(r^\lambda\right)$} \cdot \underline{\underline{O}}^{-1}
\label{eq:PowerOperatorR}
\end{equation}
where
\mbox{$\mbox{Diag$\left(r^\lambda\right)$} := 
\left( \delta_{i,j} r^{\lambda_j} \right)_{i,j=1}^D$}
From (\ref{eq:JurekCoordinates}) and (\ref{eq:PowerOperatorR})
we obtain following equations for the Jurek coordinates 
\mbox{$r := r_{\hat{\sigma} \lambda}$}
and \mbox{$\vec{\hat{\theta}}_r := \vec{\theta}_{ \hat{\sigma}\lambda}$}
of the vector $\hat{\sigma}\lambda$.
We denote \mbox{$\lambda_j := \Theta_j + \imath \Upsilon_j$} and we have:
\begin{equation}
\lambda = \left( \sum\limits_{j=1}^D 
\frac{\left|\tilde{\sigma}_j\right|^2}{r^{2 \Theta_j}} \right)^{-\frac{1}{2}}
=
\left| r^{-\underline{\underline{E}}^T} \hat{\sigma} \right|^{-1}
\quad\mbox{and}\quad
\vec{\hat{\theta}}_r = \lambda 
\left( 
 \sum\limits_{j=1}^D \underline{\underline{O}}_{i,j} \frac{\tilde{\sigma}_j}{r^{\lambda_j}}
\right)_{i=1}^D
=
\lambda r^{-\underline{\underline{E}}^T} \hat{\sigma}
\label{eq:JurekCoordinatesGenerCase}
\end{equation}
where \mbox{$\tilde{\sigma}_i := \underline{\underline{O}}^{-1}_{i,j} \hat{\sigma}_j$}
are projections of the unit vector $\hat{\sigma}$ onto the eigenvectors
of the stable index (rows of the matrix $\underline{\underline{O}}^{-1}$
or columns of the matrix $\underline{\underline{O}}$).
If the unit vector is proportional to the $l$th eigenvector then
$\tilde{\sigma}_i = \delta_{i,l}$ and from  (\ref{eq:JurekCoordinatesGenerCase})
we get $r = \lambda^{\Theta_l^{-1}}$ and 
\mbox{$\vec{\hat{\theta}}_r = \lambda^{-\imath \Upsilon_l/\Theta_l} \hat{\sigma}$}.

The fractional moment reads:
\begin{eqnarray}
E\left[(\vec{\sigma}\cdot \vec{L}_t)^\beta\right] &=&
{\mathfrak C}\left(\beta,\Theta_l^{-1}\right)
\cdot
\left( \left|\vec{\sigma}\right| t^{\Theta_l}\right)^\beta
\cdot
\nonumber \\&&
\left(
\cos(\frac{\pi}{2} \beta) e^{\imath \frac{\pi}{2} \beta}
\right)
\frac{1}{2\pi}
\int\limits_0^{2\pi} d\xi \phi^{\beta \theta_l}\left( e^{-\imath \xi} \hat{\sigma} \right)
1_{
\sum\limits_{j\ne J} \tilde{\sigma}_j^2
\left|\phi\left(e^{-\imath \xi} \hat{\sigma}\right)\right|^{2\Theta_j}
\le
\left|\vec{\sigma}\right|^2 
\left|\phi\left(e^{-\imath \xi} \hat{\sigma}\right)\right|^{2\Theta_l}
}
\label{eq:FracMomGenericUnitary}
\end{eqnarray}
Here $l$ is such a number that 
\mbox{$\Theta_l - \Theta_j \ge 0$} for all $j\ne l$ such that 
$\left|\tilde{\sigma}_j\right| > 0$,
and \mbox{$J := \left\{ j \left| \Theta_j = \Theta_l \right. \right\}$},
$D^{*} = D - \mbox{card$(J)$}+1$, and 
${\mathfrak C}(\beta, \Theta_l^{-1})$ is defined in 
(\ref{eq:PreFactorFracMom}).
The moment exists for $\beta \Theta_l < 1$, meaning
if $\beta$ does not exceed the inverse biggest real part of eigenvalues of the tail index 
$\underline{\underline{E}}$ in the maximal eigenspace containing $\hat{\sigma}$.
We note the following:

\noindent[1] The moment 
is a product of a real prefactor,
the $\beta$th power of the length of the vector $\vec{\sigma}$,
a time factor $t^{\beta \Theta_l}$ and a complex prefactor.
The former prefactor is the same as in the one dimensional case
whereas the later prefactor is a complex number equal to the support
of the random variable \mbox{$\left(\vec{\sigma} \cdot \vec{L}_t\right)^\beta$}.
In particular for symmetric distributions
the later factor is real and the support of random variable
is given by \mbox{$1 + (-1)^\beta = \cos(\frac{\pi}{2} \beta) 
e^{\imath \frac{\pi}{2} \beta}$}.

\noindent[1] The integrand in (\ref{eq:FracMomGenericUnitary}) 
is related to the complement of the $J$-space, meaning
a linear span of eigenvectors whose real parts of eigenvalues
equal $\Theta_l$.

\noindent[2] The result (\ref{eq:FracMomGenericUnitary}) is in accordance with
an existence result (Theorem 8.3.10 in \cite{Meerschaert}).
However we have for the first time computed the moment in closed form
which will be useful for statistical inference for example
or for other theoretical work. 

\noindent[3] If the multiplicity of $\Theta_l$ is equal to $D$ then $\mbox{card$(J) = D$}$
and the complement of $J$ is empty
and the last term in 
(\ref{eq:FracMomGenericUnitary}) reduces to 
\begin{equation}
\left< \phi^{\beta \Theta_l} \right> 
:= (2\pi)^{-1} 
\int\limits_0^{2\pi}  d\xi 
\phi^{\beta \Theta_l} \left(e^{-\imath \xi} \hat{\sigma}\right)
\label{eq:AveragePhiSphere}
\end{equation}
because the Heaviside function in the integrand is identically equal unity.
This is like in the scaling \& rotation case.

\noindent[4] If $\phi(e^{-\imath \xi_l} \hat{\sigma}) <1$ then 
the left hand side of the equality in the subscript of the Heaviside function
is positive and the Heaviside function may not be identically equal unity.

\section{The option price \label{sec:OptionPrice}}
An option on a financial asset 
is an agreement settled at time $t$ to purchase (call) or to sell (put)
the asset at some maturity time $T$ in the future.
Here we consider European style options that can be exercised
only at maturity. This means that boundary conditions are imposed on the option 
price at maturity \mbox{$t=T$}. 
Extending the analysis to American style options that may be exercised at any time can be done by considering European style options with a different number of exercise times \cite{Dash1} 
and allowing the number of exercise times to go to infinity.

In order to minimize the risk we now divide the money available
between $N_S$ stocks $S_t$
and $N_C$ options $C(S_t;t)$.
The value of the portfolio is then:
\begin{equation}
V(t) = N_S S_t + N_C C(S_t;t)
\label{eq:portfolioDef}
\end{equation}
We may, without loss of generality, chose $N_C = 1$.

The portfolio is a stochastic process that is required to grow exponentially
with time in terms of its expectation value.
The rate of growth $r$ is the so-called "riskless" rate of interest
and is assumed to be independent of time $t$.

\subsection{Local temporal growth\label{sec:localTemp}}

Consider the distribution of deviations
\begin{equation}
{\mathfrak D}_t^{(dt)} := V(t + d t) - e^{r d t} V(t)
\label{eq:DevDefinition}
\end{equation}
between the interest \mbox{$(e^{r d t} -  1 )V(t) = (r dt + O(dt^2)) V(
t)$}
that is earned by the portfolio
and the change $V(t + d t)- V(t)$
of the price of the portfolio. 
Does a self-financing strategy exists? 
Is it possible to choose $C = C(S_t;t)$ subject to a condition
$C_T = \mbox{max$\left( S_T - K, 0 \right)$}$, for some strike price $K$,
such that the expectation
value of the deviations of the portfolio 
conditioned on the price of the stock at time $t$
equals zero? Thus we require that the deviations have no drift:
\begin{equation}
E\left[ {\mathfrak D}_t^{(dt)} \left| S_t \right. \right] = 0
\label{eq:NoDrift}
\end{equation}

In our model we assume that the above condition is satisfied only 
for an infinitesimal time change $d t$ 
and is conditioned on the value of the stock price at time 
$t$ \mbox{(local temporal growth  sec. \ref{sec:localTemp})}.
Due to limited space we are not able to include a model extension 
that assumes that the above condition is satisfied for a finite $dt$
We will present it in a future publication.

Our approach is more general than that used in financial mathematics \cite{RamaCont}
where considerations are based on the lack of arbitrage, meaning
the assumption that riskless opportunities for making money in financial transactions
do not exist.
We waive that unrealistic assumption and instead 
require the portfolio to increase exponentially with time.

%
%

From equation (\ref{eq:StockPriceFluct}) we have:
\begin{eqnarray}
S_{t + d t} - S_t &=& 
S_t \left( \exp\left[ \alpha d t + \vec{\sigma} \cdot \vec{L}_{d t} \right] -1 \right)
\label{eq:StockPriceChangeI}
\end{eqnarray}
where we have used the fact that a L\'{e}vy process is homogeneous in time,
meaning that
\begin{equation}
\vec{L}_{t + dt} - \vec{L}_t \stackrel{d}{=} \vec{L}_{d t}
\label{eq:TimeHomogeneity} 
\end{equation}
where $\stackrel{d}{=}$ in (\ref{eq:TimeHomogeneity}) means an equality is in distribution.
We note that (\ref{eq:StockPriceChangeI}) is merely a transformation
of equation  (\ref{eq:StockPriceFluct}) and not a solution to that equation.
As such equation (\ref{eq:StockPriceChangeI}) holds for infinitesimal times $dt$ only.

From (\ref{eq:StockPriceChangeI}) 
we see
that 
the expectation value of the right hand side conditioned on $S_t$
and is infinite, unless the fluctuations are Gaussian.

Therefore we will modify the log-characteristic function $\phi(\vec{\lambda})$ 
in order to ensure the finiteness of all moments.
We define:
\begin{equation}
\phi_{\epsilon}(\vec{\lambda}) := \phi(\vec{\lambda}) \exp( -\frac{\epsilon}{\left|\vec{\lambda}\right|} )
\label{eq:ProbabilitiesModifications}
\end{equation}
where $\epsilon > 0$ and replace  $\phi_\epsilon$ by $\phi$.
From now on we will work with a fictious process related to the modified log-characteristic function, 
we will solve the option pricing problem for it
and at the end of the derivation we will take the limit $\epsilon \rightarrow 0$.
After finishing the derivation we will check analytically if the result ensures a risk free portfolio.
Firstly we check that conditional expectation value of the stock price is finite. We have:

\begin{equation}
E\left[ \left. S_{t + d t} - S_t \right| S_t \right]  =
S_t
\left(
\exp\left( \alpha d t \right) 
E\left[ e^{\left(\vec{\sigma} \cdot \vec{L}_{d t}\right)} \right] -1 
\right)
=
S_t
\left(
\exp\left( \alpha d t \right) 
e^{-dt \phi(-\imath \vec{\sigma})}
-1 
\right)
\le \infty
\label{eq:ConditionalExpValue}
\end{equation}
where in the second equality in (\ref{eq:ConditionalExpValue}) we used the following identity:
\begin{eqnarray}
E\left[ e^{\left(\vec{\sigma} \cdot \vec{L}_{d t}\right)} \right]
&=&
\int\limits_{\mathbb R} dz e^z \delta\left(z - \xi\right) 
    \nu_{\vec{\sigma} \cdot \vec{L}_{dt}}(\xi) d\xi
=
\int\limits_{\mathbb R} dz e^z 
\cdot
\frac{1}{2\pi} \int\limits_{\mathbb R} dk e^{\imath k (\xi - z)}
\cdot
\nu_{\vec{\sigma} \cdot \vec{L}_{dt}}(\xi) d\xi
\label{eq:MyIdentities} \\
&=&
\int\limits_{\mathbb R} dk \delta(k + \imath)  
\tilde{\nu}_{\vec{\sigma} \cdot \vec{L}_{dt}}(k)
=
\int\limits_{-\imath + \mathbb R} dk \delta(k + \imath)  
\tilde{\nu}_{\vec{\sigma} \cdot \vec{L}_{dt}}(k)
=
e^{-dt \phi(-\imath \vec{\sigma})}
\label{eq:MyIdentities1}
\end{eqnarray}
In the first equality in (\ref{eq:MyIdentities}) we inserted a delta function into the definition
of the expectation value, in the second equality we used the integral representation of
the delta function, in the first equality in (\ref{eq:MyIdentities1}) we integrated over
$z$ and $\xi$ and we used the integral representation of the delta function
in the second equality we shifted the integration line by using the Cauchy
theorem applied to a rectangle 
\mbox{$[-R,R] \cup R+\imath [0,1] \cup -\imath + [R,R] \cup R-\imath[0,1]$}
in the limit $R \rightarrow \infty$
and in the last equality
we used (\ref{eq:FourTrNMarginalProp}) and (\ref{eq:ScalingLimitCorFourTransf}).
We make three remarks.
Firstly the delta function has been analytically continued to complex arguments, ie we have
defined it as follows:
\begin{equation}
\delta(k + \imath q) := \sum\limits_{p=0}^\infty \frac{(\imath q)^p}{p!} \delta^{(p)}\left(k\right)
\label{eq:ImaginaryDelta}
\end{equation}
Secondly we note that the result (\ref{eq:MyIdentities1}) holds only for $\epsilon > 0$ because otherwise,
all $p\ge 2$ terms in the sum (\ref{eq:ImaginaryDelta}) produce 
infinite values when integrated with  the second term in the integrand.
Thirdly we reiterate that it is the fictious, modified stock price, related to $\epsilon >0$,
 that has a finite expectation value whereas the real stock price has of course an infinite expectation value.
The option prices that we compute correspond to a fictious $\epsilon$-world where stock prices'
probabilities have been modified like in (\ref{eq:ProbabilitiesModifications})
and the real option price is obtained as a limit of the former as $\epsilon$ tends to zero
and the sequence of the fictious worlds towards our real world.
In other words the option pricing problem has indeed no solution in our real world
however it has a solution in the ``complement'' of our world by the limit of the $\epsilon$-worlds.

We will therefore construct a zero-expectation
 value
stochastic process (\ref{eq:DevDefinition}) as a linear combination (\ref
{eq:portfolioDef})
of two stochastic processes
$S_t$ and $C(S_t;t)$ that have both non-zero expectations values.
For this purpose we will analyze the probability distribution of the deviation 
variable
${\mathfrak D}_t^{(dt)}$ and work out conditions for the option price such that
the conditional expectation value 
$E\left[{\mathfrak D}_t^{(dt)}\left| S_t \right. \right]$ is equal zero.
Now we compute the deviation of the portfolio:

\begin{eqnarray}
\lefteqn{
{\mathfrak D}_t^{(dt)} = \left( V(t+dt) - V(t) \right) + V(t) \left(1 - e^{r dt}\right) }
\label{eq:DeviationResult} \\
&&
= N_S \left( S_{t+ d t} - S_t \right) +
    \left( C_{t+ d t} - C_t \right) +
V(t) \left(1 - e^{r dt}\right) 
\label{eq:DeviationResultI} \\
&&
= N_S S_t \left( e^{\alpha dt + \vec{\sigma}\cdot \vec{L}_{d t}} - 1 \right) +
    \left( C( S_{t + d t} ;t+ d t) - C_t \right) +
V(t) \left(1 - e^{r dt}\right) 
\label{eq:DeviationResultII} \\
&&
= 
N_S S_t \left( e^{\alpha dt + \vec{\sigma}\cdot \vec{L}_{d t}} - 1 \right) +
\frac{\partial C}{\partial t} d t +
\sum\limits_{m=1}^\infty \frac{1}{m!} \frac{\partial^m C}{\partial S^m}
  S_t^m \left( e^{\alpha dt + \vec{\sigma}\cdot \vec{L}_{d t}} - 1 \right)^m +
V(t) \left(1 - e^{r dt}\right) 
\label{eq:DeviationResultIII} 
\end{eqnarray} 

In (\ref{eq:DeviationResultIII}) we have expanded the price of the option
in a Taylor series to the first order in time and to all orders in the price of the option.
In that we have assumed that the price  of the option is a perfectly smooth 
function of the price of the stock. This may limit the class of solutions.
In particular, solutions may exist, where the price of the stock
is a function satisfying the H\"{o}lder condition:
\begin{equation}
 \left| C( S_{t + d t} ;t+ d t) - C_t \right| \le A \left| S_{t+dt} - S_t \right|^\Lambda
 \label{eq:HoelderCondition}
\end{equation}
for any $S_{t + dt}$ and $S_t$, a constant $A$ and a H\"{o}lder exponent $\Lambda \in [0,1)$
and thus the price of the option can be expanded in a fractional Taylor series
\cite{Samko} in powers of $S_{t + dt} - S_t$. We will seek for these solutions in future work.

The process ${\mathfrak D}_t^{(dt)}$ is a sum of infinitely many terms 
that have non-zero expectation values.
We could compute its expectation value directly using  (\ref{eq:ConditionalExpValue})
and re-sum the series.
However, we will instead 
calculate the characteristic function of the process ${\mathfrak D}_t^{(dt)}$
conditioned on the value of the process $S_t$ at time $t$.
This means that we propagate the process $S_t$ by an infinitesimal value $dt$
and we compute the characteristic function of the increment and we require
the zero value derivative of the characteristic function to be equal zero.
This technique is not new, 
see discussion about solving master equations of Markov processes in \cite{Redner}, 
and it works because
of the time-homogeneity of the process (\ref{eq:TimeHomogeneity})
and because of the fact that the parameters $\vec{\sigma}$ and $\alpha$ are constant
as a function of the process $S_t$.
The time-homogeneity follows from the infinite divisibility of the process
and thus the technique applied here also works in the generic setting of L\'{e}vy
processes. We will extend the model according to these lines in future investigations.
 
We note that ${\mathfrak D}_t^{(dt)}$ in (\ref{eq:DeviationResultIII}) is a function
of the scalar product $\vec{\sigma} \cdot \vec{L}_{d t}$ only and thus
the distribution of ${\mathfrak D}_t^{(dt)}$ is unique functional of the distribution
of the scalar product. 

We derive the distribution of ${\mathfrak D}_t^{(dt)}$ now.
We define 
\begin{equation}
D_t := \left( \partial_t C dt + V(t)\left( 1 - e^{r dt} \right)\right)
\label{eq:SomeConstantDef}
\end{equation}
we condition on the value of the fluctuation
\mbox{$\vec{\sigma} \cdot \vec{L}_{d t}$}, we use (\ref{eq:DeviationResultIII}),
and we get:
\begin{equation}
\chi_{{\mathfrak D}_t^{(dt)}\left|S_t\right.}(k) =
\exp\left( \imath k D_t \right)
\left(
1 +
\sum\limits_{m=1}^\infty
\frac{(\imath k)^m}{m!}
\sum\limits_{q=0}^\infty
{\mathfrak A}_m(q) \exp\left(-dt(\phi(-\imath q\vec{\sigma}) - q\alpha)\right)
\right)
\label{eq:CharacteristicFct}
\end{equation}
where $\phi$ is the negative logarithmic characteristic function
of the random vector $\vec{L}_1$ (see (\ref{eq:ScalingLimitCorFourTransf})).
The expectation value of the portfolio deviation
conditioned on the value of the price of the stock $S_t$ reads:
\begin{eqnarray}
\lefteqn{E\left[{\mathfrak D}_t^{(dt)}\left|S_t \right.\right] =
\left.\frac{d \chi_{{\mathfrak D}_t^{(dt)}\left|S_t\right.}(k)}{d (\imath k)}\right|_{k=0}}
\label{eq:ExpectValue0} \\
&&=
\left(
\partial_t C - r V(t) 
- (N_S + \frac{\partial C}{\partial S}) S_t 
\left( {\mathfrak S}_{1,\phi} - \alpha \right)
-
\sum\limits_{n=1}^\infty
{\mathfrak E}_n 
\frac{\partial^n C}{\partial (\log(S))^n}
\right) dt  
+ O\left(dt^2 \right) 
\label{eq:ExpectValueII} 
\end{eqnarray}
where 
\begin{equation}
{\mathfrak S}_{n,\phi} :=
\sum\limits_{q=0}^n
\left(\begin{array}{c} n \\ q \end{array} \right)
(-1)^{n - q} \phi(-\imath q \vec{\sigma})
\quad\mbox{and}\quad
{\mathfrak E}_n =
\sum\limits_{k=\mbox{max$(n,2)$}}^\infty
\frac{a_n^{(k)}}{k!}
{\mathfrak S}_{k,\phi}
\label{eq:CoeffsBlackandScholes}
\end{equation}
and the log-characteristic function $\phi(-\imath q \vec{\sigma})$ in 
(\ref{eq:ExpectValueII})
has been analytically continued to imaginary arguments.
Here the coefficients $a_n^{(k)}$ read:


\begin{equation}
a_n^{(k)} := (-1)^{k - n} (k - 1)!
\!\!\!\!\!\!\!\!\!\!\!
\sum\limits_{1 \le j_1 < \dots < j_{n-1} \le k-1}
\prod_{q=1}^{n-1} \frac{1}{j_q}
=
(-1)^{k - n}
\!\!\!\!\!\!\!\!\!\!\!
\sum\limits_{1 \le j_1 < \dots < j_{k-n} \le k-1}
\prod_{q=1}^{k-n} {j_q}
\label{eq:Coeffs}
\end{equation}
with \mbox{$a_1^{(k)} = (-1)^{k-1} (k-1)!$}.
In addition the coefficients ${\mathfrak E}_n$ satisfy:
\begin{equation}
\sum\limits_{n=1}^\infty {\mathfrak E}_n = 0
\label{eq:SumEqZero}
\end{equation}
what follows readily from the fact that 
\mbox{$\sum_{n=1}^k a^{(k)}_n = 0$}
for $k \ge 2$.

In the Gaussian case the coefficients read
\mbox{${\mathfrak S}_{n,\phi} = (-\sigma)^2\left( n \delta_{n,1} + 
n(n-1) \delta_{n,2}\right)$} and thus (\ref{eq:ExpectValueII})
yields a second order PDE.
Since the Levy distribution has been 
truncated as in (\ref{eq:ProbabilitiesModifications})
and due to (\ref{eq:EvenLogCharFct}) the result in 
(\ref{eq:CoeffsBlackandScholes}) is real.
Indeed the log-characteristic function
can be expanded in a Taylor series in even powers of the argument only
and thus its value at the negative imaginary unit is real. 
If we did not truncate we would have obtained a unrealistic complex result
as seen from (\ref{eq:JumpPdfPureScaling}).
We reiterate that the limit of truncation threshold going to zero
($\epsilon \rightarrow 0$ in (\ref{eq:ProbabilitiesModifications}))
will be taken at the end of the calculation only rather than
at intermediary stages.
If we did so at this stage we would have obtained a paradoxical
result; an infinite sum of numbers ${\mathfrak E}_n$ each of which
is infinite equals zero.

The requirement \mbox{$E\left[{\mathfrak D}_t^{(dt)}\left|S_t \right.\right] = O(dt)$}
implies a following generalized Black \& Scholes equation:

\begin{equation}
\partial_t C 
- \sum\limits_{n=2}^\infty \frac{{\mathfrak S}_{n,\phi}}{n!} 
(S_t)^n
\frac{\partial^n C}{\partial S^n} 
=
\partial_t C 
- 
\sum\limits_{n=1}^\infty
{\mathfrak E}_n 
\frac{\partial^n C}{\partial (\log(S))^n}
= r V(t)
+ (N_S + \frac{\partial C}{\partial S}) S_t 
\left( {\mathfrak S}_{1,\phi} - \alpha \right)
\label{eq:GenBlackandScholes}
\end{equation}

In order that we get further insight into the problem,
in particular in order that
we are able to solve equation (\ref{eq:GenBlackandScholes}) analytically
we find a new expression for the coefficients of the PDE
stated in the following propositions.

\noindent{\bf Proposition 1}
The coefficients 
\mbox{${\mathfrak S}_{k,\phi}$} in 
(\ref{eq:GenBlackandScholes}) read:
\begin{eqnarray}
{\mathfrak S}_{k,\phi}
&=&
\int\limits_0^\infty d\xi \tilde{\phi}(\xi) \left(-1 + e^{-\xi}\right)^k
\label{eq:FracPowerCoeffs1} 
\label{eq:FracPowerCoeffs3} 
\end{eqnarray}
for \mbox{$m \in {\mathbb N}$}.
Here $\tilde{\phi}(\xi)$ is the inverse Laplace transform
of the log-characteristic function of $\vec{L}_1$ or the L\'{e}vy 
measure of the process $\vec{L}_t$. We have:
\begin{equation}
\tilde{\phi}(\xi) := \frac{1}{2\pi \imath}
\int\limits_{\imath {\mathbb R}} dz
e^{\xi z} \phi(-\imath z \vec{\sigma})
\quad\mbox{,}\quad
\phi(-\imath z) = \int\limits_{{\mathbb R}_+} d\xi e^{-\xi z} \tilde{\phi}(\xi)
\label{eq:InverseLaplaceTransform}
\end{equation}
In the pure scaling case for $D=1$ the inverse Laplace transform $\tilde{\phi}$ reads:
\begin{equation}
\tilde{\phi}(\vec{\xi}) = \sigma^{D \mu}
\frac{1}{2\pi \imath}
\int\limits_{\imath {\mathbb R}} dz
e^{\xi z} \left|z\right|^{D \mu}
= \frac{\sigma^{D \mu}}{2 \cos(D \mu \pi/2)}
\left(I_{+,\xi}^{-D\mu} + I_{-,\xi}^{-D \mu}\right)\left[ \delta  \right](\xi)
=
\frac{\sigma^{D\mu}}{\Gamma(-D\mu) 2 \cos(\frac{\pi}{2} D\mu)}
\frac{1}{\xi^{D\mu + 1}}
\label{eq:LogCharFctPureScal}
\end{equation}
where \mbox{$I_{\pm,x}^{-\mu} = {\mathcal D}_{\pm,x}^{\mu}$}
and the later operators are Marchaud whole axis fractional derivatives 
One has to bear in mind that since the function $\phi(-\imath z \vec{\sigma})$
may be in general unbounded as a function of $z$ and thus the quantity 
$\tilde{\phi}(\vec{\xi})$ is in general not a function
but a functional. 

\noindent{\bf Proposition 2}
The coefficients 
\mbox{${\mathfrak E}_{n}$} in 
(\ref{eq:GenBlackandScholes}) read:
\begin{eqnarray}
{\mathfrak E}_n = 
\left\{
\begin{array}{cc}
(-1) 
\int\limits_{{\mathbb R}_+} d\xi \tilde{\phi}(\xi)
\left(e^{-\xi} - 1 + \xi\right)
&\quad\mbox{if $n=1$} \\
(-1)^n 
\int\limits_{{\mathbb R}_+} d\xi \tilde{\phi}(\xi)
\frac{\xi^n}{n!}
&\quad\mbox{if $n>1$}
\end{array}
\right.
\label{eq:MyCoeffsFinal}
\end{eqnarray}
From (\ref{eq:MyCoeffsFinal}) and (\ref{eq:LogCharFctPureScal})
we see that the coefficients are infinite if $D\mu < 2$.

We proceed as follows to solve the PDE (\ref{eq:GenBlackandScholes}).
In the definition (\ref{eq:MyCoeffsFinal}) 
of the coefficients 
${\mathfrak E}_n$ we truncate the upper limit of integration at some threshold
value then we solve the generalized Black\& Scholes equation 
(\ref{eq:GenBlackandScholes}) analytically by Fourier transforming with respect
to $\log(S)$ and at the end we take the limit of the truncation threshold 
to infinity. Note that this step is essential.
Indeed, as seen from (\ref{eq:InverseLaplaceTransform}) and
from (\ref{eq:ProbabilitiesModifications})
it is not clear if the inverse Laplace transform $\tilde{\phi}$
related to  the truncated Levy distribution
diminishes fast enough away from the origin and thus
if the integral in  (\ref{eq:MyCoeffsFinal}) exists.
We accomplish this task in section (\ref{sec:FinalRes}).
Prior to doing that we describe how we will compute the number stocks as follows.

We define a utility function $U$ of the portfolio as
a functional of the price of the stock viz:
\begin{equation}
U := \int\limits_0^T V(\xi) d\xi
= \int\limits_0^T \left( N_S S_\xi + C(S;\xi) \right) d\xi
\label{eq:UtilityFunction} 
\end{equation}
and require (\ref{eq:UtilityFunction}) to be minimal.
We do not investigate here the mathematical subtleties 
concerned with the existence of the stochastic integral (\ref{eq:UtilityFunction}).
The necessary condition is that the variation $\delta U$
with respect to the price of the stock functional is zero. We have:
\begin{equation}
\delta U := 
\int\limits_0^T \left( N_S + \frac{\partial C(S;\xi)}{\partial S} \right) \delta S  d\xi
=0
\label{eq:VariationUtilityFct}
\end{equation}
what yields that 
\begin{equation}
N_S = - \frac{\partial C(S;\xi)}{\partial S}
\label{eq:NumberofStocks}
\end{equation}
as in the Gaussian case.
We note that this choice of the number of stocks ensures the self-financing property
of the portfolio.
Indeed in the Cox-Ross-Rubinstein binary tree model
in discrete time 
one considers a portfolio composed of a stock and a bond
and one derives the number of stocks
by requiring contingent claim replication,
meaning an equality of the portfolio and the claim 
with probability one (see \cite{CRRMusielaRutkowski} for example).
The later result is essentially the same as that in (\ref{eq:NumberofStocks}).

\noindent{\bf Comments} 
We have derived a PDE for the option price that ensures that the derivative 
of the expectation value of the portfolio with compounded interest is zero
\begin{equation}
\mbox{lim$_{d t \rightarrow 0}$} 
\frac{E\left[V(t + dt) - e^{r dt} V(t)\left| S_t \right.\right]}{d t} = 0
\end{equation}
without making any assumptions about the relationship between
the drift of the stock price $\alpha$ and the riskless rate of interest $r$.
We differ in that from standard models in 
financial mathematics \cite{ContOption,HurstRachev},
models that assume at the outset that $\alpha = r$. 

\subsection{Final result\label{sec:FinalRes}} 
We solve the generalized Black\& Scholes 
equation analytically. 
Inserting (\ref{eq:NumberofStocks}) 
into the second equality in (\ref{eq:GenBlackandScholes})
we get:
\begin{equation}
\partial_t C = r C - r \frac{\partial C}{\partial x}
+ 
\sum\limits_{n=1}^\infty
{\mathfrak E}_n 
\frac{\partial^n C}{\partial x^n}
\label{eq:GenBlackandScholesAnalyticSol0}
\end{equation}
where $x = \log(S_t)$.
The coefficients ${\mathfrak E}_n$ are defined in (\ref{eq:MyCoeffsFinal})
with the upper limit of integration being truncated at some threshold value.
Since the coefficients do not depend on $x$ the PDE 
(\ref{eq:GenBlackandScholesAnalyticSol0}) is converted 
into a Ordinary Differential Equation (ODE) 
by taking a Fourier transform of the option price with respect to $x$. This gives:
\begin{equation}
\partial_t \tilde{C}(k;t) = 
\left(r + H(k)\right)
\tilde{C}(k;t)
\label{eq:GenBlackandScholesAnalyticSol0FourierTrafo}
\end{equation}
where 
\begin{equation}
C(x;t) := (2\pi)^{-1} \int\limits_{{\mathbb R}} dx \tilde{C}(k;t) e^{-\imath k x}
\label{eq:FourierTrafo}
\end{equation}
and
\begin{equation}
H(k) := 
\left[
  r \imath k 
+
\sum\limits_{n=1}^\infty
{\mathfrak E}_n
(-\imath k)^{n}
\right]
= r \imath k 
 + {\mathcal V}(k)
\label{eq:Hamiltonian}
\end{equation}
We insert (\ref{eq:MyCoeffsFinal}) into (\ref{eq:Hamiltonian})
and obtain the following expression for the function ${\mathcal V}(k)$
that we call `the Hamiltonian` after Hagen Kleinert \cite{PathIntegralsOptions}.
We have:
\begin{equation}
{\mathcal V}(k) = 
\int\limits_{{\mathbb R}_+} d\xi \tilde{\phi}(\xi)
\left[
(-1 + e^{-\xi}) (\imath k) + e^{\imath k \xi} - 1
\right] =
(-1)^2
\int\limits_{{\mathbb R}_+} d\xi I^{(2)}_{-,\xi} \left[\tilde{\phi}\right](\xi)
\left[
(\imath k e^{-\xi})  - k^2 e^{\imath k \xi}
\right]
\label{eq:HamiltonianForm1}
\end{equation}
We see that the integrals in (\ref{eq:HamiltonianForm1})
exist. Therefore the limit of the truncation threshold in these integrals going 
to infinity can be performed at this stage. This is what we do now
and assume hereafter the whole positive real axis in the integration
in (\ref{eq:HamiltonianForm1}).
From (\ref{eq:InverseLaplaceTransform}) we obtain the Hamiltonian 
\begin{equation}
{\mathcal V}(k) = \left(\imath k\right) \phi(-\imath \vec{\sigma}) + \phi(-k\vec{\sigma})
\label{eq:HamiltonianPureScal}
\end{equation}
Now we come back to equation (\ref{eq:GenBlackandScholesAnalyticSol0FourierTrafo})
which we solve subject to an initial condition at maturity as follows:
\begin{equation}
\tilde{C}(k;t) = \tilde{C}(k;T)
\exp\left\{-\left(r + H(k)\right)
\tau\right\}
\label{eq:GenBlackandScholesAnalyticSol0FourierTrafoI}
\end{equation}
where $\tilde{C}(k;T)$ is the Fourier transform
of the option payoff $C(x;T)$ at maturity $T$ and $\tau := T-t$
is the time to maturity.
This payoff reads:
\begin{equation}
C(x; T) = \left\{
\begin{array}{cc}
\mbox{max$(e^x - K, 0)$} & \mbox{for a call} \\
\mbox{max$(K - e^x, 0)$} & \mbox{for a put}
\end{array}
\right.
\label{eq:MaturityCondit} 
\end{equation}
where $K$ is the strike price.
The Fourier transform of the payoff is easily computed and it reads:
\begin{eqnarray}
\tilde{C}(k;T) := \int\limits_{\mathbb R} d\xi C(\xi;T) e^{\imath k \xi}
=
\left\{
\begin{array}{cc}
K^{\imath k + 1} \left(-\frac{1}{\imath k + 1}  + \frac{1}{\imath k} + 2\pi \delta(k-\imath) \right) & \mbox{for a call} \\
K^{\imath k + 1} \left(+\frac{1}{\imath k + 1}  - \frac{1}{\imath k} - 2\pi \delta(k-\imath) \right) & \mbox{for a put} 
\end{array}
\right.
\label{eq:FourierTrafoOptionMaturityPayoff} 
\end{eqnarray}
We insert (\ref{eq:FourierTrafoOptionMaturityPayoff})
into (\ref{eq:GenBlackandScholesAnalyticSol0FourierTrafoI}) and invert the
Fourier transform for a call. We have: 
\begin{eqnarray}
\lefteqn{
C(x;t) := \frac{1}{2\pi} \int\limits_{{\mathbb R}} dk \tilde{C}(k;t) e^{-\imath k x}}
\label{eq:InverseFourierTrafoOption0} \\
&=&
S_t \!\!\!\!\!
\int\limits_{-\infty}^{-\log(m) + r \tau} \!\!\!\!\!\!d\xi
e^{-\xi}
\nu_{\vec{\sigma}\cdot\vec{L}_\tau}
(\xi + \phi(-\imath \vec{\sigma}) \tau)
-
K e^{-r \tau} 
\int\limits_{-\infty}^{-\log(m) + r \tau} \!\!\!\!\!\!d\xi
\nu_{\vec{\sigma}\cdot\vec{L}_\tau}
 (\xi + \phi(-\imath \vec{\sigma}) \tau)
%
\label{eq:InverseFourierTrafoOption6} \\
&=&
S_t N_1(d_1)
-
K e^{-r \tau} 
N_2(d_1)
\label{eq:InverseFourierTrafoOption7}
\end{eqnarray}
We recall that here $\nu_{\vec{\sigma}\cdot\vec{L}_\tau}(\xi)$
is the probability density function of the fluctuation term 
$\vec{\sigma}\cdot\vec{L}_\tau$.

In (\ref{eq:InverseFourierTrafoOption7}) we have changed the integration variables
and simplified the result. Here we defined:
\begin{equation}
N_1(d) := \int\limits_{-\infty}^{d} d\xi
e^{-\xi} e^{\tau \phi(-\imath \vec{\sigma})}
\nu_{\vec{\sigma}\cdot\vec{L}_\tau}(\xi)
\;\mbox{,}\;
N_2(d) := \int\limits_{-\infty}^{d} d\xi 
\nu_{\vec{\sigma}\cdot\vec{L}_\tau}(\xi)
\;\mbox{and}\;
d_1 = -\log(m) + \tau \left( r + \phi(-\imath \vec{\sigma}) \right)
\label{eq:ConstDefs}
\end{equation}
In the limit $D\mu \rightarrow 2_{-}$ the density 
$\nu_{\vec{\sigma}\cdot\vec{L}_\tau}(\xi)$
in (\ref{eq:InverseFourierTrafoOption6})
goes into a Gaussian with mean zero and variance $2 \sigma \tau$ and
(\ref{eq:InverseFourierTrafoOption7})
goes into the Gaussian Black\& Scholes equation,
see e.g. eqs. (1.6),(1.7) on page 8 in \cite{RamaCont}.
We end this section by stating the price of the portfolio. We have: 
\begin{equation}
V(t) = N_S S_t  + C(x;t) =
- K e^{-r \tau} N_2(d_1)
\label{eq:PortfolioCheck4}
\end{equation}
Since the last factor on the right-hand side in (\ref{eq:PortfolioCheck4})
depends implicitly on $S_t$ the unconditional expectation value
of the portfolio does not increase exponentially as required.
Therefore the solution (\ref{eq:InverseFourierTrafoOption7}) is only an approximation.
However since the Gaussian Black \& Scholes equation is a particular case of
(\ref{eq:InverseFourierTrafoOption7}) it turns out that it 
is also only an approximation.

The factors in (\ref{eq:ConstDefs}) are complex which is of course unrealistic.
The reason for that is the following.
In our approach we assumed that the time change $dt$ is infinitesimally
small rather asumming it to be finite at the outset and taking the limit 
$dt \rightarrow 0$  at the end of the calculation.
We have checked that the later procedure leads to a real result
which has the same form as in (\ref{eq:InverseFourierTrafoOption7})
except that the L\'{e}vy density 
\mbox{$\nu_{\vec{\sigma}\cdot\vec{L}_\tau}(\xi)$}
goes into an inverse Fourier transform of 
\mbox{$\exp\left( 1/2 \sigma^\mu \Gamma(\imath k + \mu)/\Gamma(\imath k) \right)$}
evaluated at \mbox{$\xi - r \tau \phi(-\imath \sigma)$}
which essentially amounts to replacing the expression
$-k^\mu + \imath k$ by $\Gamma(\imath k + \mu)/\Gamma(\imath k)$
in some intermediate calculations.
The inverse Fourier transform in question is essentially
equal to the 
$\log(S_\tau) - \alpha\tau$ process 
probability density function 
evaluated at the argument \mbox{$\xi - r \tau \phi(-\imath \sigma)$}.
Therefore the price of the option is a discounted present value 
of the maturity payoff under a risk-neutral probability measure
where the measure in question is related to the compensated log-price
process $log(S_t) - \alpha t$.
Thus we have proven that the risk-neutral option pricing method
holds in the generic setting of  operator stable processes.

Expressions (\ref{eq:ConstDefs}) are difficult to deal with in numerical calculations.
Indeed the ``typical width'' of the Fourier transform of the 
L\'{e}vy density is $\sigma^\mu \tau$. Since this quantity is small,
meaning of the order of $10^{-2}$
for stock daily data and for times to maturity of the order of hundreds
of days, the use of ``primitive'' methods like Romberg quadratures
for evaluating the Fourier integrals requires a very high precision
of calculation that is much bigger than the precision of the estimated
parameters. Therefore we propose to use a more sophisticated method
for the numerical integrations. This method is described in the Appendix.

\section{Conclusions\label{sec:Concl}}
We have applied the technique of characteristic functions
to the problem of pricing an option on a stock that is driven
by operator stable fluctuations. 
We have developed a technique to ensure that the expectation value of the portfolio 
grows exponentially with time.
In doing this we have not, unlike other authors, 
made any assumptions about the analytic properties of the log-characteristic function 
of the stock price process. 
Instead we have expressed all results in terms of the 
characteristic function of the operator \mbox{stable fluctuation $\vec{L}_1$}.

 Subsequent to successful numerical tests, we ought then to be able to price analytically not only European options but also exotic options with a finite number of different exercise times. This should also allow us to price American style options by allowing the number of exercise times to become infinite.

We may also compute the $99$th percentile of the probability 
distribution of the deviation of the portfolio (Value at Risk) 
as a function of $\vec{\sigma}$ and of the log-characteristic function $\phi$ of the random vector $\vec{L}_1$.
The Value at Risk will be expressed as an integral equation involving the conditional characteristic function of the portfolio deviation (\ref{eq:CharacteristicFct}). The resulting integrals will be carried out by means of the Cauchy complex integration theorem. 

The results of these calculations will be reported in a future publication.

\section{Acknowledgments}
This work resulted from research conducted within the SFI Basic Research Grant 04/BR/0251.
We are grateful to Mark Meerschaert, Stefan Thurner, Christoli Bieli and Krzysztof Urbanowicz
for useful discussions.

\section{Appendix } We explain how the integrals from the Levy density in 
(\ref{eq:InverseFourierTrafoOption7}) are computed
numerically in the case $D=1$. 
Note that in this case equation (\ref{eq:InverseFourierTrafoOption7}) 
can be written as follows:
\begin{equation}
C(x;t) = S_t N_{\sigma L_\tau}^{(1)}(d;z) 
       - K e^{-r \tau} N_{\sigma L_\tau}^{(0)}(d;z)
\label{eq:GenBlackScholesFinalSimplif}
\end{equation}
where \mbox{$z := z_r + \imath z_i = \phi(-\imath \sigma) \tau$} is a complex number,
$d := - \log(m) + r \tau$ is a real number and
\begin{equation}
N_{\sigma L_\tau}^{(s)}(d;z) :=
\int\limits_{-\infty}^d d\xi 
e^{-s \xi}
\nu_{\sigma L_\tau} (\xi + z) 
\label{eq:FactorsSimplif}
\end{equation}
for $s=0,1$.
We note that the factors (\ref{eq:FactorsSimplif}) have a following
integral representation that lends itself to numerical computations in a 
straightforward manner. We have:
\begin{eqnarray}
\lefteqn{
N_{\sigma L_\tau}^{(s)}(d;z) =
\int\limits_{-\infty}^d d\xi 
e^{-s \xi}
\cdot
\frac{1}{2\pi}
\int\limits_{\mathbb R} dk e^{-\imath k (\xi + z)} 
\cdot
e^{-\tau \phi(\sigma k) }
\mathop{=}_{A\rightarrow \infty}
\frac{1}{2\pi}
\int\limits_{\mathbb R} dk e^{-\imath k z} 
\cdot
\int\limits_{-A}^d d\xi 
e^{-(s + \imath k) \xi}
e^{-\tau \phi(\sigma k) } }
\label{eq:FactorProof1} \\
&&
\!\!\!\!\!\!\!\!\!\!\!\!\!\!
\mathop{=}_{A\rightarrow \infty}
\frac{1}{2\pi}
\int\limits_{\mathbb R} dk e^{-\imath k z} 
\cdot
\left(
\frac{-e^{-\imath \theta d} + e^{\imath \theta A}}{\imath \theta}
\right)
\cdot
e^{-\tau \phi(\sigma k) }
=
\mathop{=}_{A\rightarrow \infty}
\frac{1}{2\pi}
\int\limits_{{\mathbb R} + \imath s} dk e^{-\imath k z} 
\cdot
\left(
\frac{-e^{-\imath \theta d} + e^{\imath \theta A}}{\imath \theta}
\right)
\cdot
e^{-\tau \phi(\sigma k) }
\label{eq:FactorProof2} \\
&&
\!\!\!\!\!\!\!\!\!\!\!\!\!\!
\mathop{=}_{A\rightarrow \infty}
\frac{1}{2\pi}
\int\limits_{\mathbb R} d\theta 
e^{-\imath (\theta + \imath s) z}
\cdot
\left(
\frac{-e^{-\imath \theta d} + e^{\imath \theta A}}{\imath \theta}
\right)
\cdot
e^{-\tau \phi(\sigma(\theta + \imath s)) }
\label{eq:FactorProof3} \\
&&
\!\!\!\!\!\!\!\!\!\!\!\!\!\!
=
\frac{e^{s z}}{2}
\left[
e^{-\tau \phi(\imath s\sigma)}
+
\frac{1}{\pi}
\int\limits_0^\infty d\theta e^{\theta z_i}
\left[
\frac{\sin(\theta (d+z_r))}{\theta} 
  {\mathcal M}_1(\theta)
+ 
\imath \frac{\cos(\theta (d+z_r))}{\theta} 
  {\mathcal M}_2(\theta)
\right]
\right]
\label{eq:FactorProof4}
\end{eqnarray}
In the first equality in (\ref{eq:FactorProof1}) we expressed the L\'{e}vy
stable density through its Fourier transform and in the second equality in  
(\ref{eq:FactorProof1}) we changed to order of integration.
In the first equality in (\ref{eq:FactorProof2}) we integrated over $\xi$
and we defined 
$\theta = k - \imath s$
and in the second equality in (\ref{eq:FactorProof2}) we shifted the integration line
by $\imath s$ in the complex plane. In doing this we used the Cauchy theorem
applied to a rectangular contour composed of an interval $\left[-R,R\right]$,
of that interval shifted by $\imath s$ and of sections perpendicular to the real axis
that complete the contour. In the limit $R\rightarrow \infty$ the integrals over the
later sections vanish.
In (\ref{eq:FactorProof3}) we factorize-d the integral
and in (\ref{eq:FactorProof4}) we performed the limit $A\rightarrow \infty$
by using the identity:
\begin{equation}
\mbox{lim$_{A\rightarrow \infty}$} 
\frac{e^{\imath \theta A}}{\imath \theta} = \pi \delta(\theta)
\quad \mbox{for} \quad \theta \in {\mathbb R}
\end{equation}
where
\begin{equation}
{\mathcal M}_1(\theta) := 
 \sum\limits_{p=\pm 1}   e^{-\tau \phi(\sigma(p \theta + \imath s))} =
 \quad\mbox{and}\quad
{\mathcal M}_2(\theta) := 
  \sum\limits_{p=\pm 1} p e^{-\tau \phi(\sigma(p \theta + \imath s))}
\label{eq:MFactorsDef}
\end{equation}
In the pure scaling case in one dimension, from (\ref{eq:PureScalnubeta}),
we have \mbox{$\phi(k) = \phi_\pm k^\mu$} and thus:
\begin{equation}
{\mathcal M}_1(\theta) := 
 \sum\limits_{p=\pm 1}   e^{- \phi_\pm \tau (\sigma l)^\mu \cos(\mu \phi_p) }
 \left( \cos(\alpha_p) - \imath \sin(\alpha_p) \right)
 \quad\mbox{and}\quad
{\mathcal M}_2(\theta) := 
 \sum\limits_{p=\pm 1}   p e^{-\phi_\pm \tau (\sigma l)^\mu \cos(\mu \phi_p) }
 \left( \cos(\alpha_p) - \imath \sin(\alpha_p) \right)
\end{equation}
where
\begin{equation}
l := \sqrt{ \theta^2 + s^2 },
\cos(\phi_p) = \frac{p \theta}{l}, 
\sin(\phi_p) = \frac{s}{l},
 \quad\mbox{and}\quad
 \alpha_p = \phi_\pm \tau (\sigma l)^\mu \sin(\mu \phi_p) 
\end{equation}
Since, as seen from (\ref{eq:MFactorsDef}),
\mbox{${\mathcal M}_1(0) = 2 e^{-\tau \phi(\imath \sigma s)}$}
and
\mbox{${\mathcal M}_2(0) = 0$}
the integral in (\ref{eq:FactorProof4}) is clearly finite
the result can be used for numerical calculations.
In the Gaussian case $\mu = 2$ we have 
\mbox{$\alpha_p = 2 \phi_\pm \tau \sigma^2 s p \theta$} and thus
\begin{equation}
{\mathcal M}_1(\theta) := 
2 \cos(2 \phi_\pm s \sigma^2 \theta \tau) e^{-\phi_\pm \tau \sigma^2 \left( \theta^2 - s^2 \right)}
\quad\mbox{and}\quad
{\mathcal M}_2(\theta) := 
-2 \imath \sin(2 \phi_\pm s \sigma^2 \theta \tau) e^{-\phi_\pm \tau \sigma^2 \left( \theta^2 - s^2 \right)}
\label{eq:MFactorsDefGaussian}
\end{equation}
and \mbox{$z = \tau \phi(-\imath \sigma) = z_r + \imath z_i = -\tau \phi_\pm \sigma^2$}. 
Inserting (\ref{eq:MFactorsDefGaussian}) into (\ref{eq:FactorProof4}) gives:
\begin{eqnarray}
N_{\sigma L_\tau}^{(s)}(d;z) &=&
\frac{1}{2} + \frac{1}{2 \pi}
\int\limits_{\mathbb R} d\theta e^{-\tau \phi_\pm ( \sigma \theta )^2}
\cdot
\frac{\sin(\theta  e )}{\theta}
=
\frac{1}{2} + \frac{1}{2 \pi}
\int\limits_{\mathbb R} d\theta e^{-\tau \phi_\pm ( \sigma \theta )^2}
\cdot
\left(
\frac{1}{2} \int\limits_{-e}^e d\eta e^{-\imath \eta \theta} 
\right)
\label{eq:NFactorGaussianCheck} \\
&=&
\frac{1}{2} +
\frac{1}{2} \int\limits_{-e}^e d\eta 
\nu_{\sigma L_\tau}(\eta) 
=
\int\limits_{-\infty}^e d\eta
\nu_{\sigma L_\tau}(\eta) 
=
\int\limits_{-\infty}^{\frac{e}{\sqrt{2\phi_\pm \sigma^2 \tau}}} d\eta  \frac{1}{\sqrt{2\pi}} \exp\left\{-\frac{1}{2} \eta^2 \right\}
\label{eq:NFactorGaussianCheck1}
\end{eqnarray}
where \mbox{$e = -\log(m) + \tau \left(r + \phi_\pm \sigma^2(2 s - 1)\right)$}. 
From (\ref{eq:NFactorGaussianCheck1}) we see that the factors coincide with those in the Gaussian Black\& Scholes formula.

\end{document}